\documentclass[10pt, a4paper, english, twoside]{article}
\usepackage[utf8]{inputenc}
\usepackage{ifpdf}
\usepackage{anysize}
\usepackage[colorlinks=true, linkcolor=blue, urlcolor=black]{hyperref}
\usepackage{enumerate}
\usepackage{amssymb, amsmath, amsthm}
\usepackage{mathtools}
\usepackage[x11names]{xcolor, colortbl}
\usepackage[T1]{fontenc}
\usepackage{listings}
\usepackage{subfig}
\ifpdf
  \usepackage[pdftex]{graphicx}
  \DeclareGraphicsExtensions{.pdf,.png,.jpg}
\else
  \usepackage[dvips]{graphicx}
  \DeclareGraphicsExtensions{.eps}
  \fi  
\definecolor{roig}{rgb}{1, 0.5, 0.5}
\definecolor{verd}{rgb}{0.7, 1, 0.5}
\definecolor{groc}{rgb}{1, 1, 0.4}
\definecolor{blau}{rgb}{0.2, 0.5, 1}
\definecolor{cyan}{rgb}{0, 1, 1}
\definecolor{ombra}{rgb}{0.95, 0.95, 1}
\definecolor{gris}{gray}{0.85}
\newcolumntype{g}{>{\columncolor{gris}}c}
\newif\ifdepurar
\depurarfalse
\ifdepurar
\usepackage[notref,notcite]{showkeys}

\newcommand{\blau}[1]{{\color{blau}#1}}
\else

\newcommand{\blau}[1]{#1}
\fi
\newcommand{\map}[2]{\varphi_{#1}^{[#2]}}
\newcommand{\opA}{\mathcal{L}_A}
\newcommand{\opB}{\mathcal{L}_B}
\title{Efficient symplectic integrators for cubic and quartic potentials}
\author{Alejandro Escorihuela-Tom\`as \\[2ex]
  {\small\it Departament de Matem\`atiques and IMAC, Universitat Jaume I, 12071-Castell\'o, Spain}\\{
\small\it email: alescori@uji.es}}
\makeindex

\begin{document}
\maketitle
\begin{abstract}
  We present a set of new, efficient high-order symplectic methods designed for Hamiltonian systems with cubic or quartic potentials.
  By demonstrating that polynomial potentials require fewer order conditions,
  we develop schemes that outperform both standard symmetric compositions of second-order methods and existing RKN splitting methods.
  Numerical results confirm their improved efficiency over state-of-the-art alternatives found in the literature.
\end{abstract}
\ifdepurar
\tableofcontents
\fi
\section{Introduction}
We are interested in constructing symplectic integrators for Hamiltonian systems of the form:
\begin{equation}\label{eq:HT}
  H(q,p) = T(p) + V(q), \quad \text{with} \quad T(p) = \frac{1}{2m}\sum_{i=1}^d p_i^2,
\end{equation}
where $q_i, p_i, m \in \mathbb{R}$ for $i = 1, \dots, d$, with $d$ denoting the dimension of the system, $q=(q_1,\dots, q_d)$, $p=(p_1,\dots,p_d)$ and $V(q)$ is a polynomial potential.
For Hamiltonian systems, the equations of motion are given by:
\begin{equation*}
  \dot{q}_i = \frac{\partial H}{\partial p_i}, \qquad \dot{p}_i = -\frac{\partial H}{\partial q_i},
\end{equation*}
which, in the case of a quadratic kinetic energy as in \eqref{eq:HT}, can be rewritten as:
\begin{equation}\label{eq:qdotdot}
  \ddot{q}_i = -\frac{1}{m} \frac{\partial V}{\partial q_i}.
\end{equation}
This second-order system corresponds to an ordinary differential equation of the form:
\begin{equation}\label{eq:rkneq}
  \ddot{y} = g(y),
\end{equation}
where, if $V(q)$ is a polynomial of degree $n$, the function $g(q)$ is a polynomial of degree $n - 1$. For equations of these form, specially tailored numerical methods can be constructed.
In this work, we focus on such equations with the additional assumption that $g(q)$ is a polynomial of degree up to 3, which corresponds to quartic potentials.
Furthermore, we require our methods to preserve the geometric structure of the original problem;
in particular, we search symplectic integrators.

The defining feature that distinguishes symplectic integrators from other numerical methods is that each integration step corresponds to the application
of a symplectic transformation. As a consequence, such integrators preserve phase-space volume and ensure that the energy error remains bounded.
In particular, even over very long integration times, no systematic growth of the energy error is observed \cite{blanes2024smf}.

In the next section, we introduce splitting methods, which will be used to construct our symplectic schemes.
In Section 3, we derive the order conditions that the coefficients of a splitting method must satisfy in order to achieve a prescribed order of accuracy.
Section 4 presents the newly developed integration methods, while in Section 5 we compare their performance with that of existing symplectic
integrators reported in the literature. Finally, Section 6 provides some concluding remarks.
\section{Splitting methods}
One efficient approach to constructing symplectic integrators for Hamiltonian problems is through the use of splitting methods \cite{blanes2025aci}. We begin by reviewing this technique, followed by a section on the derivation of order conditions. We then present new integrators and provide numerical experiments to compare their performance against the most efficient schemes currently available in the literature.

Splitting methods are applicable when the system $\dot{x} = f(x)$ can be decomposed into two or more components, each of which can be integrated exactly. In our setting, we consider a two-part splitting:
\begin{equation}\label{eq:split}
  \dot{x} = f(x) = f_A(x) + f_B(x), \quad \text{with} \quad x(t_0) = x_0,
\end{equation}
where the flows corresponding to each subproblem $\dot{x} = f_i(x)$ admit exact solutions denoted by $\map{t}{i}(x_0) = x(t_0 + t)$ for $i = A, B$. By composing these flows with appropriately chosen coefficients, we can construct an integrator for the full problem:
\begin{equation*}
  \psi_t = \map{a_0t}{A} \circ \map{b_0t}{B} \circ \map{a_1t}{A} \circ \map{b_1t}{B} \circ \cdots \circ \map{a_{s-1}t}{A} \circ \map{b_{s-1}t}{B},
\end{equation*}
where $s$ denotes the number of stages, with each stage corresponding to a evaluation of $\map{t}{A} $ and $\map{t}{B}$. Denoting by $\varphi_t(x_0)$ the exact solution, the integrator is said to be of order $r$ if:
\begin{equation*}
  \psi_t(x_0) = \varphi_t(x_0) + \mathcal{O}(t^{r+1}).
\end{equation*}

It is well known that compositions of symplectic maps yield symplectic. Since the flows arising in Hamiltonian splitting are themselves symplectic, the overall integrator inherits this property. Splitting schemes preserve the qualitative features of the dynamical system and exhibit favorable long-term energy behavior \cite{hairer2006gni}. One of the most widely used second-order splitting schemes is the Strang splitting:
\begin{equation}\label{eq:strang}
  \mathcal{S}_t^{[2]} = \map{t/2}{A} \circ \map{t}{B} \circ \map{t/2}{A}.
\end{equation}

In this study, we focus on splitting methods with palindromic sequence of coefficients:
\begin{equation}\label{eq:splitsim}
  \psi_t = \map{a_0t}{A} \circ \map{b_0t}{B} \circ \map{a_1t}{A} \circ \map{b_1t}{B} \circ \cdots \circ \map{a_{s-2}t}{A} \circ \map{b_{s-2}t}{B} \circ \map{a_{s-1}t}{A}, \qquad \text{with } a_j = a_{s-j-1},\ b_j = b_{s-j-2},
\end{equation}
since such schemes are symmetric in time, i.e., $\psi_{-t} \circ \psi_t \equiv \mathrm{Id}$ \cite{blanes2025aci}, where $\mathrm{Id}$ denotes the identity transformation.
\blau{In other words, advancing by a step $t$ and then by a step $-t$ returns us to the initial point, thereby preserving, by construction, one of the properties of the exact flow.}

Problems of the form $\ddot{y} = g(y)$ (of which our system is a particular case) can be treated via splitting by introducing an auxiliary variable $v = \dot{y}$ and rewriting the system as:
\begin{equation}\label{eq:splitgy}
  \dot{x} \equiv (\dot{y}, \dot{v})^\top = (v, 0)^\top + (0, g(y))^\top = f_A(y, v) + f_B(y, v),
\end{equation}
where the flows admit exact solutions:
\begin{equation*}
  \map{t}{A}(y, v) = (y + t v, v)^\top, \qquad \map{t}{B}(y, v) = (y, v + t g(y))^\top.
\end{equation*}

  As mentioned above, for this kind of problems ($\ddot{y}=g(y)$), specially adapted numerical methods can be constructed.
  For instance, Runge--Kutta methods tailored to such equations are known as Runge--Kutta--Nyström (RKN) methods \cite{hairer1993sod}; by analogy, splitting methods adapted to this class of problems are commonly referred to as RKN splitting methods \cite{blanes2002psp}.
Splitting schemes of up to sixth order, and up to eighth order in the RKN framework, have been constructed (see, e.g., \cite{blanes2024smf}).
However, one of the main challenges in high-order methods is the exponential growth in the number of order conditions required \cite{mclachlan2002sm}.
To address this, a highly effective alternative is to construct high-order symplectic integrators using symmetric composition techniques
based on a symmetric second-order method, such as the Strang splitting \eqref{eq:strang}:
\begin{equation}\label{eq:compss}
  \mathcal{SS}_t = \mathcal{S}_{\alpha_0 t}^{[2]} \circ \mathcal{S}_{\alpha_1 t}^{[2]} \circ \cdots \circ \mathcal{S}_{\alpha_{s-1} t}^{[2]}, \qquad \text{with } \alpha_j = \alpha_{s - j}.
\end{equation}
This strategy has led to the development of efficient high-order symplectic integrators \cite{sofroniou2005dos}, which we use as benchmarks to evaluate the performance of the new methods proposed in this article.

Having introduced the class of problems under consideration and the general framework of splitting methods, we now proceed to derive the order conditions required to construct high-order schemes.
\section{Order conditions}
Splitting methods for problems of the form $\ddot{y} = g(y)$ require fewer order conditions to achieve a given order compared to general splitting methods.
Furthermore, when $g(y)$ is a polynomial function, an additional reduction in the number of order conditions occurs.
In this work, we consider the cases where $g(y)$ is either quadratic or cubic, corresponding to cubic and quartic potentials, respectively.

As shown in \cite{blanes2025aci, blanes2024smf}, each splitting integrator can be associated with a formal series $\Psi(t)$ of differential operators.
In the case of the time-symmetric splitting integrator introduced in the previous section (see equation~\eqref{eq:splitsim}), we have:
 \begin{equation*}
   \Psi(t) = e^{a_{s-1}t\opA}e^{b_{s-2}t\opB}\cdots e^{b_{0}t\opB}e^{a_{0}t\opA},
 \end{equation*}
 where $\opA$ and $\opB$ denote the Lie derivatives associated with the functions $f_A$ and $f_B$, respectively \cite{arnold1989mmo}. In the case considered here namely,
 $\ddot{y} = g(y)$ with the splitting introduced previously, these operators take the form:
 \begin{equation*}
   \opA = \sum_{i=1}^dv_i\frac{\partial}{\partial y_i},\qquad\qquad \opB = \sum_{i=1}^dg_i\frac{\partial}{\partial v_i}.
 \end{equation*}
 Associating a differential operator series $\Psi(t)$ to the numerical integrator allows us to apply the Baker--Campbell--Hausdorff (BCH) formula \cite{varadarajan1984lgl}, expressing the composition $\Psi(t)$ as the exponential of a single operator, i.e., $\Psi(t) = \exp(Z(t))$, where:
 \begin{equation}\label{eq:zt}
   Z(t)=\sum_{i \ge 1}t^i\sum_{j=1}^{n_i} \omega_{i,j}E_{i,j},
 \end{equation}
 with $E_{i,j}$ denoting a basis of independent commutators of order $i$ resulting from the application of the BCH formula, where $E_{1,1} = \opA$ and $E_{1,2} = \opB$. The integers $n_i$ represent the number of independent commutators at order $i$, and the coefficients $\omega_{i,j}$ are polynomial functions of the splitting coefficients $a_k$ and $b_k$.
 Thus, in order for the numerical integrator to be of order $r$, the operator $Z(t)$ must satisfy
 \begin{equation*}
   Z(t)=t(\opA + \opB)+\mathcal{O}(t^{r+1}),
 \end{equation*}
 which is achieved by setting $\omega_{1,1} = \omega_{1,2} = 1$ and cancelling all other functions $\omega_{i,j}$ for $i \le r$.
 These functions at low orders are of the form \cite{blanes2025aci,blanes2024smf}
 %As an example, a time-symmetric splitting integrator of the form \eqref{eq:splitsim} with order 4 must satisfy the following order conditions:
 % \begin{equation*}
 %   \sum_{i=0}^{s-1}a_i =1,\quad \sum_{i=0}^{s}b_i =1,\quad \sum_{i=0}^{s}b_i\left(\sum_{j=0}^{i}a_j\right)^2-\frac{1}{3}=0,\quad
 %   \sum_{i=0}^{s-1}a_i\left(\sum_{j=0}^{s}b_j\right)^2-\frac{1}{3}=0.
 % \end{equation*}
 \begin{gather*}
   \omega_{1,1}=\sum_{i=0}^{s-1}a_i,\quad \omega_{1,2}=\sum_{i=0}^{s}b_i =1,\quad \omega_{2,1}=\frac{1}{12}\omega_{1,1}\omega_{1,2}-\sum_{0\le i <j \le s}b_ia_j\\
   \omega_{3,1}=\frac{1}{12}\omega_{1,1}^2\omega_{1,2}-\frac{1}{2}\sum_{0\le i< j \le k \le s}a_ib_ja_k,\quad
  \omega_{3,2}=\frac{1}{12}\omega_{1,1}\omega_{1,2}^2-\frac{1}{2}\sum_{0\le i\le j < k \le s -1}b_ia_jb_k.
 \end{gather*}
 Here, the time-symmetry of the integrator implies that all coefficients $\omega_{2k,j}$ vanish automatically \cite{blanes2025aci}.
 The choice of commutators $E_{i,j}$ is not unique. Below, we list the commutators used up to order 4 in the present work:
\begin{center}
  \begin{tabular}{lllll}
    $E_{11} = \opA$;&$E_{12} = \opB$; &$E_{21} = [\opA,\opB]$;  \\%[0.3em]
    $E_{31}= [\opA,\opA,\opB]$; &$E_{32}= [\opB,\opA,\opB]$;&  \\%[0.3em]
    $E_{41}= [\opA,\opA,\opA,\opB]$; &$E_{42}= [\opB,\opA,\opA,\opB]$; &$E_{43}= [\opB,\opB,\opB,\opA]$;
  \end{tabular}
\end{center}
Here, the basic commutator is defined as $[\opA, \opB] = \opA\opB - \opB\opA$, and nested commutators are evaluated from right to left, for example, $[\opA, \opA, \opB] \equiv [\opA, [\opA, \opB]]$. For higher orders, the first element at each level is constructed as $E_{k,1} = [\opA, \dots, \opA, \opB]$ with $\opA$ appearing $k-1$ times.
Table~\ref{tau:nc} lists the number of independent commutators arising at each order $r$ for the general splitting case, denoted by $n_r$.

In the next section, we illustrate how the number of order conditions is reduced in the case $\ddot{y} = g(y)$, and how further simplification occurs when $g(y)$ is polynomial.
\subsection{Reduction of order conditions in $\ddot{y}= g(y)$}
Using the Lie derivatives introduced earlier, the commutator $[\opA, \opB]$ takes the form \cite{blanes2025aci}:
\begin{equation*}
  [\opA,\opB] = \sum_{j,i=1}^dv_j\frac{\partial g_i}{\partial y_j}\frac{\partial}{\partial v_i}-\sum_{i=1}^dg_i\frac{\partial}{\partial y_i}.
\end{equation*}
With this result, the next nested commutator can be computed as \cite{blanes2025aci}:
\begin{equation*}
  [\opB,\opA,\opB]=2\sum_{j,i=1}^dg_j\frac{\partial g_i}{\partial y_j}\frac{\partial}{\partial v_i}.
\end{equation*}
As can be seen, this commutator does not depend on the variables $v_i$.
Consequently, higher-order nested commutators such as $[\opB, \opB, \opA, \opB]$ vanish \cite{blanes2025aci,blanes2024smf},
and many of the terms appearing in expression \eqref{eq:zt} cancel automatically.
These results can be seen in Table~\ref{tau:nc}, where we list the number of independent commutators that remain at order $r$, denoted by $l_r$. In Table~\ref{tau:ns},
we show the corresponding minimal number of stages $L_r$ required by a symmetric numerical integrator when taking advantage of this structural property.
\subsection{Reduction of order conditions in $\ddot{y}= g(y)$ with polinomical $g$}
The reduction in the number of order conditions for cubic potentials was first introduced in \cite{mclachlan2003tae, iserles1998rkm}. We now derive the expressions for the nested commutators that make this simplification possible.
Using the expression for $[\opA, \opB]$, we compute:
\begin{equation*}
  [\opA,\opA,\opB]= \smashoperator{\sum_{k,j,i=1}^d}v_kv_j\frac{\partial^2 g_i}{\partial y_k \partial y_j}\frac{\partial}{\partial v_i}
  -2\smashoperator{\sum_{j,i=1}^d}v_j\frac{\partial g_i}{\partial y_j}\frac{\partial}{\partial y_i},
\end{equation*}
which corresponds to the element $E_{3,1}$ in our notation. Proceeding to compute the next element, $E_{4,1}$:
\begin{equation*}
  [\opA,\opA,\opA,\opB]= \smashoperator{\sum_{l,k,j,i=1}^d}v_lv_kv_j\frac{\partial^3 g_i}{\partial y_l\partial y_k \partial y_j}\frac{\partial}{\partial v_i}
  -3\smashoperator{\sum_{k,j,i=1}^d}v_kv_j\frac{\partial^2 g_i}{\partial y_k\partial y_j}\frac{\partial}{\partial y_i}.
\end{equation*}
From these expressions, it is straightforward to generalize the formula for the commutator $E_{n+1,1}$ as follows:
\begin{equation}\label{eq:anb}
  [\underbrace{\opA,\dots,\opA}_{n \text{ times}},\opB]=
  \smashoperator{\sum_{i_{n+1},\dots,i_1=1}^d}v_{i_{n+1}}\cdots v_{i_2}\frac{\partial^ng_{i_1}}{\partial y_{i_{n+1}}\cdots \partial y_{i_2}}\frac{\partial}{\partial v_{i_1}}
  -n \smashoperator{\sum_{i_{n},\dots,i_1=1}^d}v_{i_{n}}\cdots v_{i_2}\frac{\partial^{n-1}g_{i_1}}{\partial y_{i_{n}}\cdots \partial y_{i_2}}\frac{\partial}{\partial y_{i_1}}.
\end{equation}
From this general expression, one can immediately conclude that if $g(y)$ is a polynomial of degree $n$, say $g(y) = P_n(y)$, then:
\begin{equation}\label{eq:enpn}
  g(y) = P_n(y) \rightarrow E_{n+3,1}= [\underbrace{\opA,\dots,\opA}_{n +2\text{ times}},\opB]=0.
\end{equation}
This vanishing occurs because the $(n+1)$-th derivatives of $g(y)$ in expression \eqref{eq:anb} are identically zero.
In this work, we are interested in constructing numerical schemes for the cases $n=2,3$. We now examine these two specific cases in detail:
\begin{itemize}
\item \textbf{Cubic potential}. A cubic potential $V(q) = P_3(q)$ leads to an equation of motion $\ddot{q} = g(q)$ where $g(q) = Q_2(q)$, with $P_3$ and $Q_2$
  denoting polynomials of degree 3 and 2, respectively. According to \eqref{eq:enpn}, the first commutator that vanishes due to this condition is $E_{5,1} = [\opA, \opA, \opA, \opA, \opB] = 0$.
  This results in a reduction in the number of order conditions, in addition to that already induced by the second-order structure $\ddot{y} = g(y)$.
  In Table~\ref{tau:nc}, the number of independent surviving commutators at order $r$ is denoted by $l_r^3$, while in Table~\ref{tau:ns},
  the corresponding minimum number of stages required for a symplectic splitting integrator is indicated by $L_r^3$.
  
\item \textbf{Quartic potential}. In the case of a quartic potential $V(q) = P_4(q)$, the corresponding force function is cubic, i.e., $g(q) = Q_3(q)$.
  Thus, from \eqref{eq:enpn}, the first vanishing commutator is $E_{6,1} = [\opA, \opA, \opA, \opA, \opA, \opB] = 0$.
  The number of independent commutators at order $r$ is denoted by $l_r^4$ in Table~\ref{tau:nc}, and the minimum number of stages
  for a corresponding symplectic integrator is indicated by $L_r^4$ in Table~\ref{tau:ns}.
\end{itemize}
We now present symplectic numerical integration methods based on the splitting technique, specifically tailored for the two cases described above.

\begin{table}[!h]
   \begin{center}
    \begin{tabular}{c|ccccccccc}
      $r$        &1 &2 &3 &4 &5 &6 &7  &8  &9 \\ 
      \hline
      $s_r$      &1 &0 &1 &1 &2 &2 &4  &5  &8 \\ 
      $n_r$      &2 &1 &2 &3 &6 &9 &18 &30 &56 \\ 
      $l_r$      &2 &1 &2 &2 &4 &5 &10 &14 &25 \\     
      $l_r^3$    &2 &1 &2 &2 &3 &3 &6  &6  &10 \\
      $l_r^4$    &2 &1 &2 &2 &4 &4 &8  &10 &18 \\
    \end{tabular}
    \caption{\small{Number of independent commutators at order $r$ using different techniques: symmetric composition of second-order methods ($s_r$),
        general splitting methods ($n_r$), RKN case ($l_r$), cubic potentials ($l_r^3$), and quartic potentials ($l_r^4$).}}
    \label{tau:nc}
  \end{center}
\end{table}
\begin{table}[!h]
\begin{center}
  \begin{tabular}{c|ccccc}
    %\hline
    $r$     &2 &4 &6 &8 &10 \\ 
    \hline
    $S_r$   &1 &3 &7 &15 &31 \\
    $N_r$   &1 &3 &9 &27 &83 \\ 
    $L_r$   &1 &3 &7 &17 &42 \\ 
    $L_r^3$ &1 &3 &6 &12 &22 \\
    $L_r^4$ &1 &3 &7 &15 &33 \\
  \end{tabular}
  \caption{\small{Minimum number of stages required to achieve order $r$ in symmetric methods using different techniques: symmetric composition of second-order methods ($S_r$), general splitting methods ($N_r$), RKN case ($L_r$), cubic potentials ($L_r^3$), and quartic potentials ($L_r^4$).}}
  \label{tau:ns}
\end{center}
\end{table}

\section{New methods}
In the case of general splitting, the flows $\map{t}{A}$ and $\map{t}{B}$ are interchangeable.
However, in the specific setting considered here, the choice of splitting type \eqref{eq:splitgy}
and the resulting cancellation of certain commutators require us to distinguish between two types of compositions \cite{blanes2025aci}.
On the one hand, we consider compositions of type ABA:
\begin{equation*}
     \psi_t = \map{a_0t}{A}\circ\map{b_0t}{B}\circ\map{a_1t}{A}\circ\map{b_1t}{B}\circ\cdots\circ\map{a_{s-1}t}{A}\circ\map{b_{s-1}t}{B}\circ\map{a_{s}t}{A},\qquad\text{with }a_{j}=a_{s-j},\ b_{j}=b_{s-j-1},
 \end{equation*}
and on the other hand, compositions of type BAB:
 \begin{equation*}
     \psi_t = \map{b_0t}{B}\circ\map{a_0t}{A}\circ\map{b_1t}{B}\circ\map{a_1t}{A}\circ\cdots\circ\map{b_{s-1}t}{B}\circ\map{a_{s-1}t}{A}\circ\map{b_{s}t}{B},\qquad\text{with }a_{j}=a_{s-j-1},\ b_{j}=b_{s-j}.
 \end{equation*}
The number of stages is counted using the FSAL (First Same As Last) property \cite{blanes2002psp, blanes2022rkn}, which allows the last flow of one integration step to be reused as the first flow of the next.

We denote the newly constructed methods using a two-letter code. The first letter is either $\mathcal{C}$ or $\mathcal{Q}$, indicating whether the method is designed for cubic or quartic potentials, respectively. The second letter specifies the composition type: $\mathcal{A}$ for ABA and $\mathcal{B}$ for BAB compositions.

% The coefficients of the new methods presented in this work are listed in Tables~\ref{tau:metC} and~\ref{tau:metQ}.
We have selected those methods that showed the best efficiency; additional methods and simulation data are available at the following link:
\begin{center}
  https://doi.org/10.5281/zenodo.19772477
\end{center}
Interestingly, the most efficient methods we obtained were of ABA type.

To solve the order conditions, we employed Python libraries \cite{python3rm}, specifically \texttt{numpy} and \texttt{scipy} \cite{numpy,scipy}. The latter provides wrappers for MINPACK routines \cite{minpack}, which implement Powell's hybrid algorithm \cite{press1992nri}. Among the solutions obtained, we selected those with the smallest $\ell_1$ norm of the coefficients, defined as:
\begin{equation*}
  \ell_1 = \sum_{i=0}^{s} |a_i| + \sum_{i=0}^{s-1} |b_i| \qquad \text{(ABA)} \qquad \text{or} \qquad \ell_1 = \sum_{i=0}^{s-1} |a_i| + \sum_{i=0}^{s} |b_i| \qquad \text{(BAB)}.
\end{equation*}
In cases where more stages than the theoretical minimum (see Table~\ref{tau:ns}) were used,
the additional degrees of freedom were exploited to minimize the $\ell_1$ norm of the coefficients using
an arc pseudo-continuation technique as described in \cite{alberdi2019aab}.
\begin{itemize}
\item \textbf{Cubic potential}.  
In this case, the complete exploration at order 4 has already been carried out in \cite{iserles1998rkm}.  
For this reason, we begin our search at order 6, exploring both ABA and BAB compositions with stage counts ranging from 6 to 11.  
With 6 stages, we recover the method presented in \cite{iserles1998rkm}, known as the 6-stage SRKNQ scheme.  
The most efficient method we found corresponds to an 11-stage composition of ABA type.  
At orders 8 and 10, the best-performing integrators were also of ABA type and achieved with the minimum number of stages required.

\item \textbf{Quartic potential}.  
For quartic potentials, we begin at order 8, where the benefit of reduced stage count becomes noticeable.  
We explore ABA and BAB compositions from 15 to 21 stages, and the most efficient result is obtained with a 19-stage ABA composition.  
At order 10, several solutions were found, but none of them outperformed the $\mathcal{SS}$ methods of order 10 in terms of efficiency.
\end{itemize}
\begin{table}[!h]
  \centering
  \renewcommand\arraystretch{1.1}
  \begin{tabular}{ll}
    \multicolumn{2}{c}{$\mathcal{CA}_{11}^{[6]}$}\\
    \hline
    $a_0= 0.042694933980191700$ &\qquad $b_0= 0.162759370295069398$ \\
    $a_1= -0.037632511090066230$ &\qquad $b_1= -0.0319763989469851006$ \\
    $a_2= 0.219222765218418850$ &\qquad $b_2= 0.193876531329159266$ \\
    $a_3= 0.172304948988733900$ &\qquad $b_3= 0.153628407800734632$ \\
    $a_4= -0.866412097755661100$ &\qquad $b_4= -0.000976972075829582387$ \\
    $a_5= \frac{1}{2}-\sum_{i=0}^4a_i$&\qquad $b_5 = 1-2\sum_{i=0}^4b_i$\\
    %% $a_5= 0.9698219606583829$ &\qquad $b_5= 0.045378123185901176$ \\      
    \hline
    \multicolumn{2}{r}{$\ell_1=5.748$}\\
    & \\
    \multicolumn{2}{c}{$\mathcal{CA}_{12}^{[8]}$}\\
    \hline
    $a_0= 0.249757865893252399$ &\qquad $b_0= 0.213843067589222296$ \\
    $a_1= 0.00573645298706788704$ &\qquad $b_1= -0.19020545357715192$ \\
    $a_2= -0.205722262874388455$ &\qquad $b_2= 0.152112874905099611$ \\
    $a_3= 0.205575388768639098$ &\qquad $b_3= 0.230725630134443253$ \\
    $a_4= -0.271817217898894439$ &\qquad $b_4= -0.0190391037211012732$ \\
    $a_5= 0.489641667780589551$ &\qquad $b_5= \frac{1}{2}-\sum_{i=0}^4b_i$ \\
    $a_6= 1-2\sum_{i=0}^5a_i$ &\qquad \\   
    % $a_5= 0.489641668225369$ &\qquad $b_5= 0.11256298465493847$ \\
    % $a_6= 0.0536562106751578$ &\qquad \\
    \hline
    \multicolumn{2}{r}{$\ell_1=4.747$}\\
    & \\
    \multicolumn{2}{c}{$\mathcal{CA}_{22}^{[10]}$}\\
    \hline
    $a_0= 0.0449093524320847274$ &\qquad $b_0= 0.12962054755858581$ \\
    $a_1= 0.206215756187292034$ &\qquad $b_1= 0.28538465968498975$ \\
    $a_2= 0.283122139467291742$ &\qquad $b_2= 0.447546675818306977$ \\
    $a_3= -0.0467135692211750245$ &\qquad $b_3= -0.223172488593771899$ \\
    $a_4= 0.191753353873687573$ &\qquad $b_4= -0.168587763117298738$ \\
    $a_5= -0.30261984583857599$ &\qquad $b_5= -0.152330297495996805$ \\
    $a_6= 0.403507309497133283$ &\qquad $b_6= 0.0357528743187445468$ \\
    $a_7= -0.411843791732104965$ &\qquad $b_7= 0.109308117709693938$ \\
    $a_8= 0.532623166827988789$ &\qquad $b_8= -0.0278888594069430645$ \\
    $a_{9}= -0.672573105536616394$ &\qquad $b_9= 0.0874141139842583174$ \\
    $a_{10}= 0.585892787986557426$ &\qquad $b_{10}= \frac{1}{2}-\sum_{i=0}^9a_i$ \\
    $a_{11}= 1-2\sum_{i=0}^{10}a_i$ &\qquad \\      
    % $a_{10}= 0.5858927879860629$ &\qquad $b_{10}= -0.02304758046060429$ \\
    % $a_{11}= -0.6285471078864806$ &\qquad \\   
    \hline
    \multicolumn{2}{r}{$\ell_1=11.372$}
  \end{tabular}
  \caption{\small{Symmetric methods of the form (\ref{eq:splitsim}) of order six, eight, and ten for cubic potentials, along with the corresponding $\ell_1$ norms of the coefficients.}}
  \label{tau:metC}
\end{table}

\begin{table}[!h]
  \centering
  \renewcommand\arraystretch{1.1}
  \begin{tabular}{ll}
    \multicolumn{2}{c}{$\mathcal{QA}_{19}^{[8]}$}\\
    \hline   
    $a_0= 0.017198824867539785$ &\qquad $b_0=0.0575638169530652679$ \\
    $a_1= 0.102745265073641400$ &\qquad $b_1=0.107137729341092432$\\
    $a_2= 0.184885965868561040$ &\qquad $b_2=0.0185107295436068784$\\
    $a_3=-0.219582702854432000$ &\qquad $b_3=-0.0195107639666207008$\\
    $a_4= 0.114544474826412145$ &\qquad $b_4=-0.0766782391416861259$\\
    $a_5=-0.0380031092293327718$ &\qquad $b_5=0.134395138847794067$\\
    $a_6= 0.112642307911035366$ &\qquad $b_6=0.110620560134710358$\\
    $a_7= 0.121499361307126452$ &\qquad $b_7=0.137518059034892416$\\
    $a_8= 0.169053217437104066$ &\qquad $b_8=-0.0368744345203615847$\\
    $a_9= \frac{1}{2}-\sum_{i=0}^8a_i$&\qquad $b_9 = 1-2\sum_{i=0}^8b_i$\\
    \hline
    \multicolumn{2}{r}{$\ell_1=3.822$}
  \end{tabular}
  \caption{\small{Coefficients and $\ell_1$ norm of the symmetric method of the form (\ref{eq:splitsim}) of order eight for quartic potentials.}}
   \label{tau:metQ}
\end{table}
\section{Numerical examples}
Below is the list of methods with which we will compare the methods we have obtained:
\begin{itemize}
\item $\mathcal{NA}_{14}^{[6]}$: Splitting RKN of order 6 with 14 stages obtained by Blanes \& Moan \cite{blanes2002psp}.
\item $\mathcal{NB}_{18}^{[8]}$: Splitting RKN of order 8 with 18 stages obtained by Blanes et al. \cite{blanes2022rkn}.
\item $\mathcal{SS}_{19}^{[8]}$ and $\mathcal{SS}_{35}^{[10]}$: $\mathcal{SS}$ methods of order 8 with 19 stages and of order 10 with 35 stages respectively obtained by Sofroniou \& Spaletta \cite{sofroniou2005dos}.
\end{itemize}
We will also compare with $\mathcal{CA}_{6}^{[6]}$, an existing 6th-order method with 6 stages for problems of this type \cite{iserles1998rkm}.

We now present two test problems involving cubic potentials and two involving quartic potentials.
\subsection{Hénon--Heiles}
One of the most studied problems with cubic potentials is the Hénon--Heiles system \cite{henon1964tao}.
\begin{equation*}%\label{eq:hehe}
  H(q,p)=\frac{1}{2}(p_1^2+p_2^2)+\frac{1}{2}(q_1^2+q_2^2) + q_1^2q_2 - \frac{q_2^3}{3},
\end{equation*}
with $q_i,p_i\in\mathbb{R}$.
% Un problema que ha sigut considerat, entre altres, com a model per caracteritzar la transició al caos Hamiltonià.
% Nosaltres resoldrem aquest problema amb les condicions inicials:
A problem that has been considered, among others, as a model to characterize the transition to Hamiltonian chaos.
We will solve this problem with the initial conditions
\begin{equation*}
  % q_1(0) = \frac{\alpha}{2},\ q_2(0)=0,\ p_1(0)=0,\ p_2(0)=\frac{\alpha}{4},\qquad\text{with }0\le\alpha\le 1.
  q_1(0) = \frac{1}{4},\ q_2(0)=0,\ p_1(0)=0,\ p_2(0)=\frac{1}{8},
\end{equation*}
conditions that ensure that the system lies within the integrable regime of the problem.

\begin{figure}[!h]
  \begin{center}
    \subfloat[]{
      \label{fig:1a}
      \includegraphics[width=8cm]{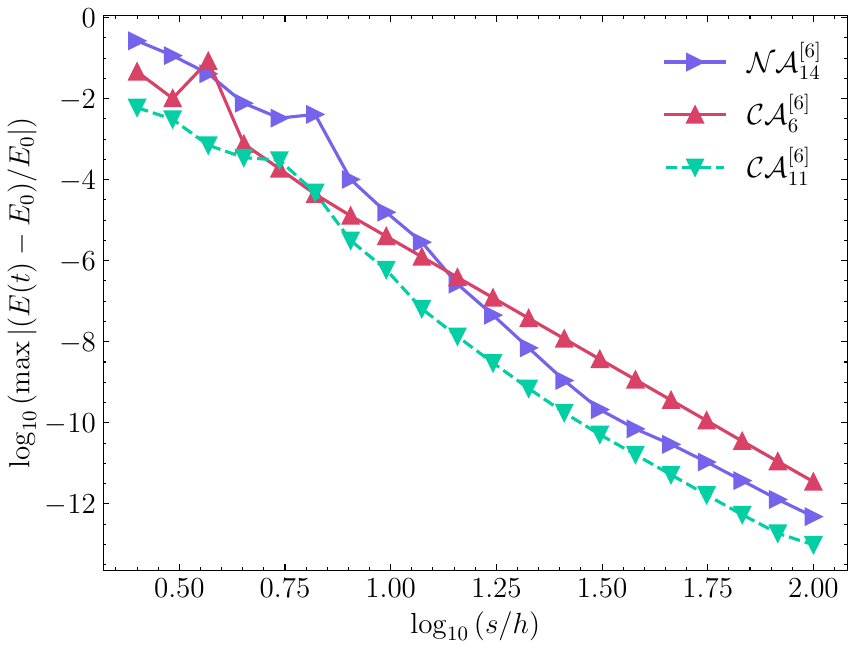}}
    \subfloat[]{
      \label{fig:1b}
      \includegraphics[width=8cm]{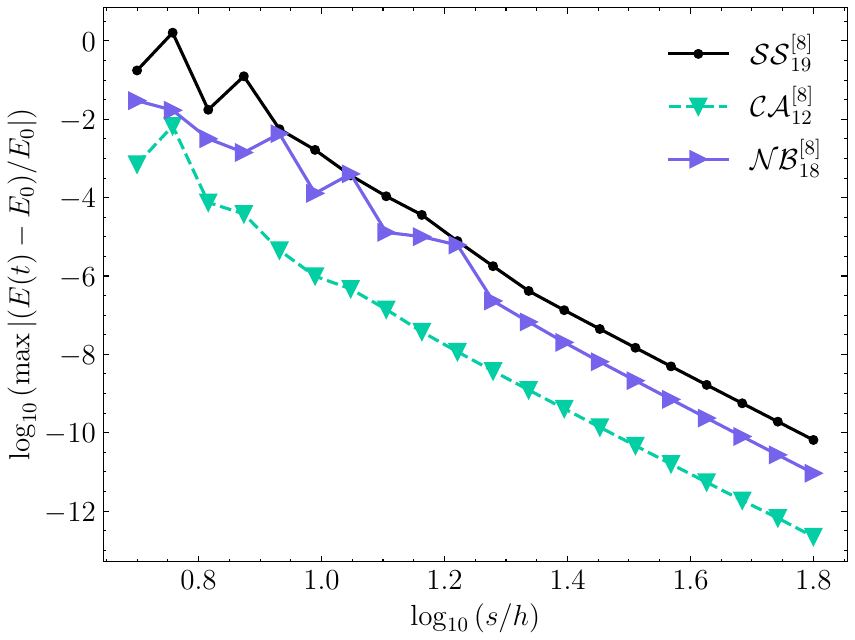}}
    \\ 
    \subfloat[]{
      \label{fig:1c}
      \includegraphics[width=8cm]{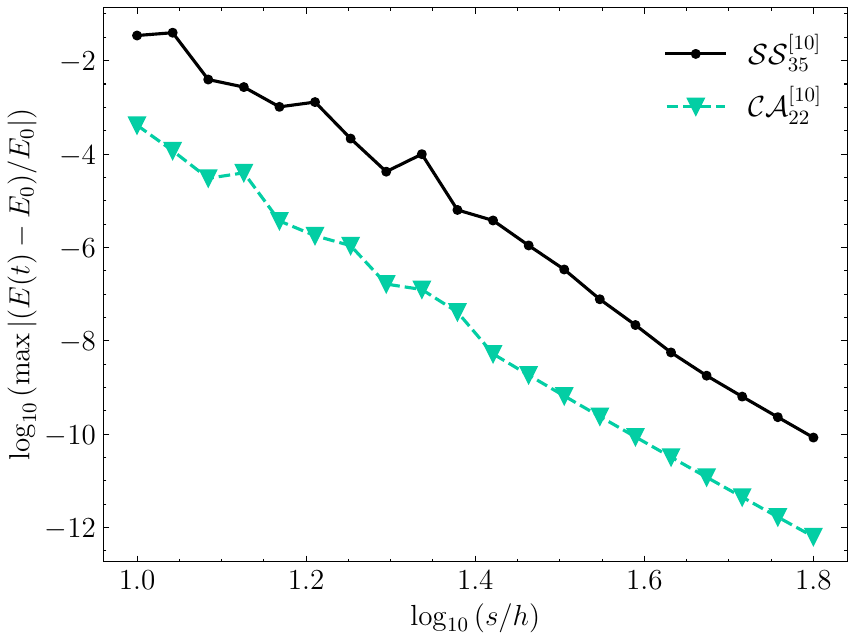}}
    \caption{Efficiency diagrams for the Hénon--Heiles problem with initial conditions $(q_1, q_2, p_1, p_2) = (1/4, 0, 0, 1/8)$ and final time $t_f = 1000$. \textbf{(a)} 6th-order methods; \textbf{(b)} 8th-order methods; and \textbf{(c)} 10th-order methods.}
    \label{fig:1}
  \end{center}
\end{figure}
Figure~\ref{fig:1} displays efficiency diagrams for the newly developed methods at orders 6, 8, and 10, respectively.  
All simulations were performed with a final integration time of $t_f = 1000$, and comparisons were made against the most efficient existing symplectic methods for each order, as well as against $\mathcal{CA}_{6}^{[6]}$.
\begin{figure}[!h]
  \begin{center}
    \subfloat[]{
      \label{fig:1a_eci}
      \includegraphics[width=8cm]{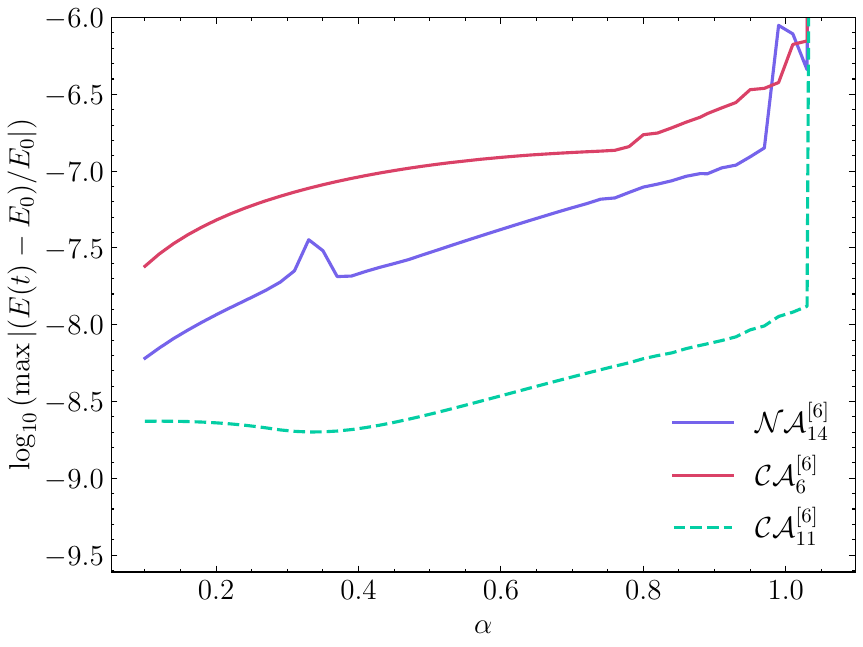}}
    \subfloat[]{
      \label{fig:1b_eci}
      \includegraphics[width=8cm]{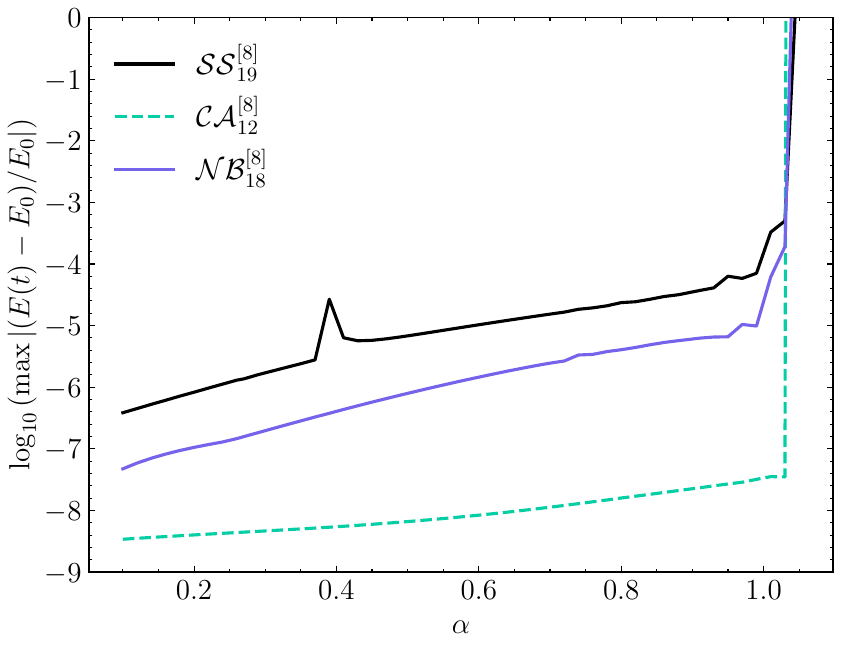}}
    \\ 
    \subfloat[]{
      \label{fig:1c_eci}
      \includegraphics[width=8cm]{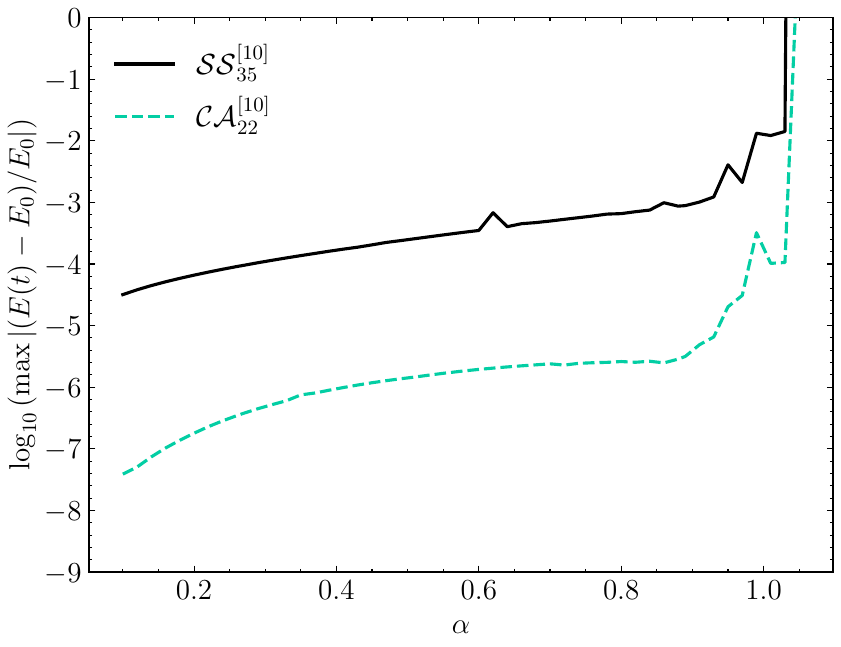}}
    \caption{Energy conservation error as a function of the initial conditions for final time $t_f = 1000$.
      All methods have been implemented using the same computational cost with $s/h=10^{5/4}$.
     \textbf{(a)} 6th-order methods; \textbf{(b)} 8th-order methods; and \textbf{(c)} 10th-order methods.}
    \label{fig:1_eci}
  \end{center}
\end{figure}

This problem is commonly used to study the transition to chaos. For energy values $E_0<\frac{1}{12}$, the system is integrable; $\frac{1}{12}<E_0<\frac{1}{6}$, chaotic behavior begins to emerge; and for $E_0>\frac{1}{6}$, the system becomes fully chaotic \cite{henon1964tao,hairer2006gni}. To explore different dynamical regimes, we parametrize the initial conditions as follows:
\begin{equation*}
  q_1(0) = \frac{\alpha}{2},\ q_2(0)=0,\ p_1(0)=0,\ p_2(0)=\frac{\alpha}{4},
\end{equation*}
so that the choice $\alpha = \frac{1}{2}$ recovers the initial conditions considered previously. In this parametrization,
the two thresholds delimiting the regimes described above are given by $\alpha \approx 0.7303$ and
$\alpha \approx 1.0328$ respectively.
Figure~\ref{fig:1_eci}  illustrates the energy conservation error for different methods as a function of the parameter $\alpha$.
The results show that the newly proposed methods are more efficient across the entire range of $\alpha$ for which the system remains integrable.
All simulations have been performed with the same computational cost, namely $s/h=10^{5/4}$ and with a final integration time of $t_f=1000$.
The results show that the newly proposed integrators consistently outperform the currently known symplectic methods in terms of efficiency.

\begin{figure}[!h]
  \begin{center}
    \includegraphics[width=8cm]{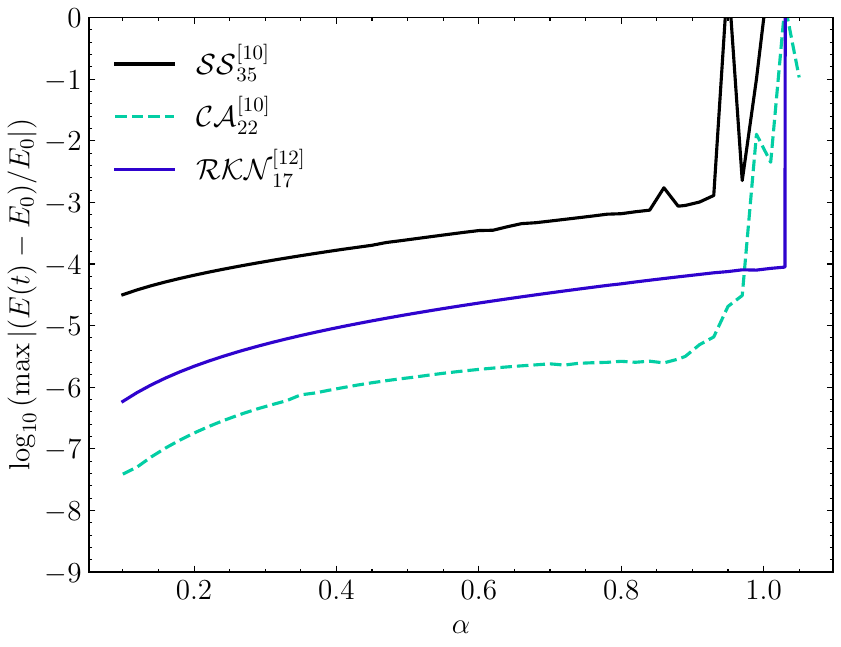}
    \caption{Same as in Figure~\ref{fig:1c_eci}, but with $t_f=10^7$ and including the non-symplectic method $\mathcal{RKN}_{17}^{[12]}$.}
    \label{fig:ecirkn}
  \end{center}
\end{figure}

We can go one step further and perform a comparison over long integration times with a non-symplectic
method.
A method that yields very good performance is the $\mathcal{RKN}_{17}^{[12]}$ scheme,
which is of order 12 and consists of 17 stages. This is an explicit Runge--Kutta method
of Nyström type (RKN),
originally introduced in \cite{dormand1987hoe}, with coefficients taken from \cite{rkn12web}.
In Figure~\ref{fig:ecirkn}, we present three simulations comparing the two tenth-order symplectic methods with the twelfth-order RKN method, for a final integration time $t_f=10^7$ and $s/h=10^{5/4}$.
The results show that the new method proposed in this work achieves higher accuracy.

\subsection{Random cubic potential}
We now present an experiment involving a randomly generated cubic potential:
\begin{equation*}
  H(q,p)=\frac{1}{2}\sum_{i=1}^dp_i^2 +  \sum_{i,j=1}^d\mu_{ij} q_iq_j + \sum_{i,j,k=1}^d\nu_{ijk} q_iq_jq_k,
\end{equation*}
where $p_i, q_i, \mu_{ij}, \nu_{ijk} \in \mathbb{R}$. The coefficients $\mu_{ij}$ (for $i \ne j$) and $\nu_{ijk}$
are randomly generated from a uniform distribution in the interval $[-\frac{1}{2}, \frac{1}{2}]$, and $d$ denotes the system dimension.
The system is constructed as a perturbation of the harmonic oscillator, with $\mu_{ii} = \frac{1}{2}$ for all $i$.
Initial conditions are also drawn from a uniform distribution in the interval $[0, \frac{1}{5}]$.

\begin{figure}[!h]
  \begin{center}
    \subfloat[]{
      \label{fig:2a}
      \includegraphics[width=8cm]{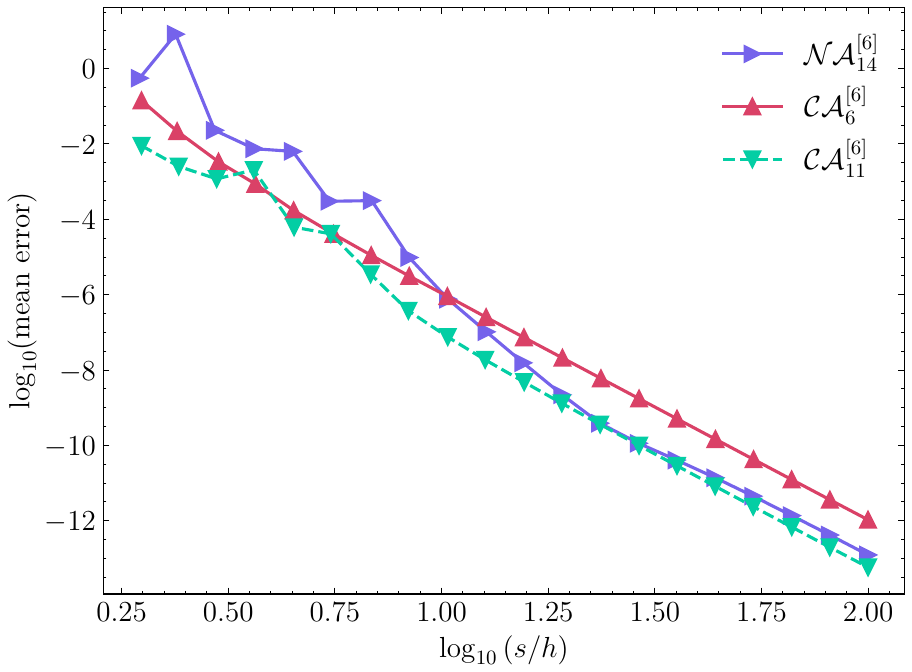}}
    \subfloat[]{
      \label{fig:2b}
      \includegraphics[width=8cm]{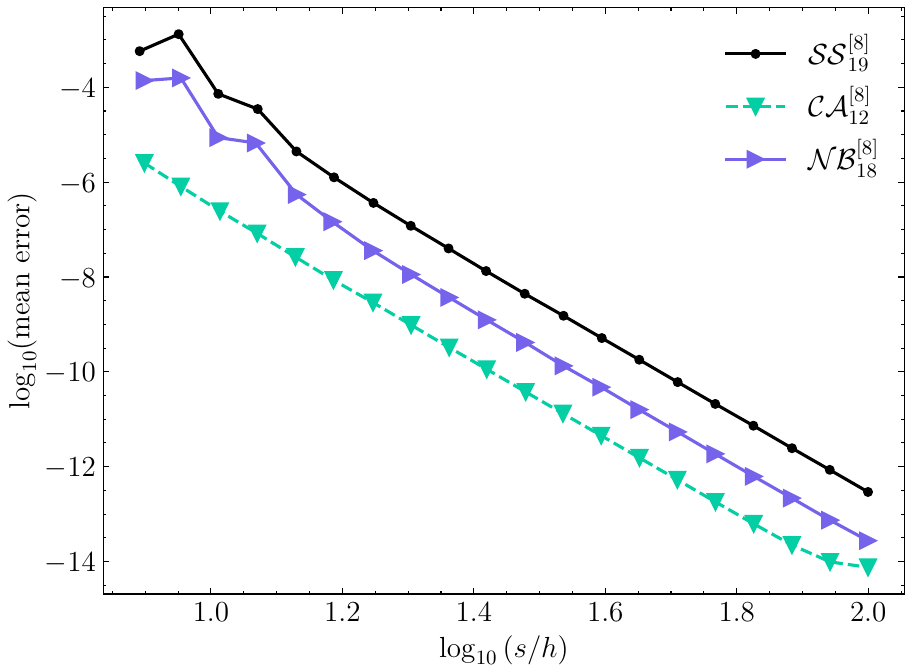}}
    \\ 
    \subfloat[]{
      \label{fig:2c}
      \includegraphics[width=8cm]{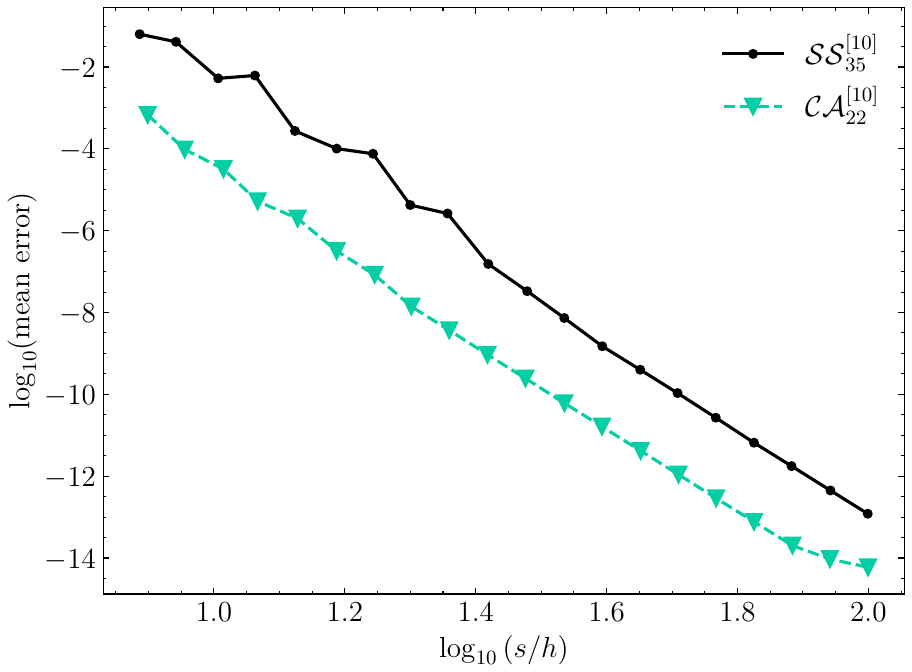}}
    \caption{Efficiency diagrams for the random cubic potential with randomly generated initial conditions, $d=10$, and final time $t_f = 100$.
      The result shown is the average error over 10 simulations.
      \textbf{(a)} 6th-order methods; \textbf{(b)} 8th-order methods; and \textbf{(c)} 10th-order methods.}
    \label{fig:2}
  \end{center}
\end{figure}

To generate Figure~\ref{fig:2}, we performed 10 simulations using different random realizations of $\mu_{ij}$, $\nu_{ijk}$,
and initial conditions.
The average behavior was computed across all simulations using the following metric:
\begin{equation}\label{eq:mean_error}
  \text{mean error :}\quad \frac{1}{M}\sum_{m=1}^{M}\log_{10}\left(\text{max}\left|\frac{E_m(t)-E_m(0)}{E_m(0)}\right|\right),
\end{equation}
where $M$ is the number of simulations.

\subsection{Random quartic potential}
As in the cubic potential case, we can generate a quartic potential as a perturbation of the harmonic oscillator:
  \begin{equation*}
    H(q,p)=\frac{1}{2}\sum_{i=1}^dp_i^2 + \sum_{i,j=1}^d\mu_{ij} q_iq_j + \sum_{i,j,k=1}^d\nu_{ijk} q_iq_jq_k+ \sum_{i,j,k=1}^d\rho_{ijkl} q_iq_jq_kq_l,
  \end{equation*}
  where $p_i, q_i, \mu_{ij}, \nu_{ijk}, \rho_{ijkl} \in \mathbb{R}$.
  The diagonal terms are set to $\mu_{ii} = \frac{1}{2}$, while $\mu_{ij}$ for $i \ne j$, $\nu_{ijk}$, and $\rho_{ijkl}$
  are pseudo-random numbers uniformly distributed in the interval $[-\frac{1}{2}, \frac{1}{2}]$.

\begin{figure}[!h]
  \begin{center}
    \subfloat[]{
      \label{fig:4a}
      \includegraphics[width=8cm]{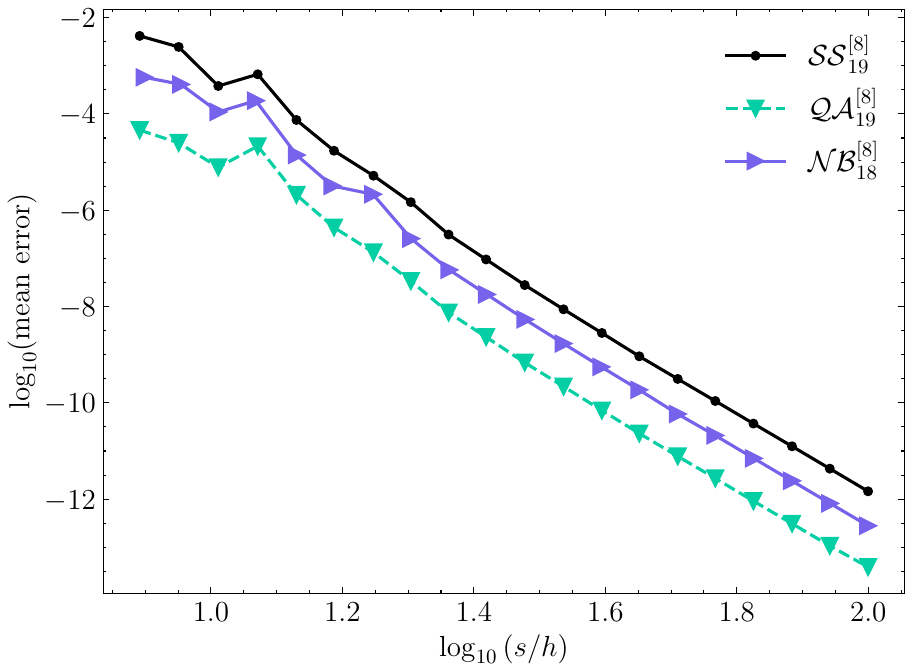}}
    \subfloat[]{
      \label{fig:4b}
      \includegraphics[width=8cm]{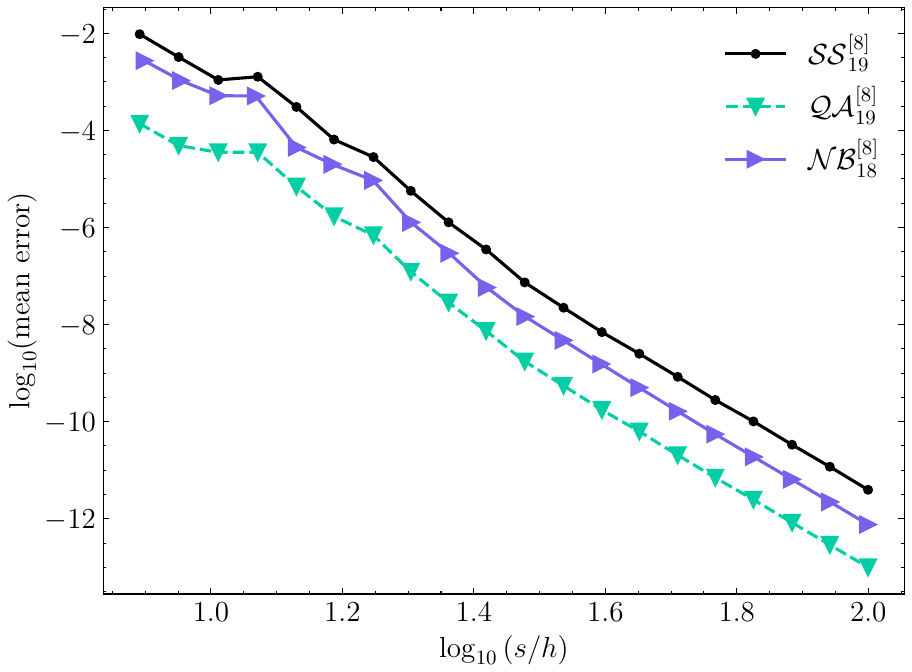}}
    \caption{Efficiency diagrams for the random quartic potential with randomly generated initial conditions and final time $t_f = 100$ with 8th-order methods.
      The result shown is the average error over 10 simulations.
      \textbf{(a)} Dimension $d = 5$; \textbf{(b)} Dimension $d = 10$.}  
    \label{fig:4}
  \end{center}
\end{figure}

Figure~\ref{fig:4} shows an efficiency diagram for our 8th-order method, compared with the best-performing $\mathcal{SS}$-type method and the most efficient 8th-order RKN method for this class of problems.

As in the cubic case, initial conditions are sampled uniformly from the interval $[0, \frac{1}{5}]$,
and the mean error is computed using formula~\eqref{eq:mean_error}, averaged over 10 simulations.
We consider two different system sizes: $d=5$ in Figure~\ref{fig:4a}, and $d=10$ in Figure~\ref{fig:4b}.

\subsection{Schrödinger equation with quartic potential}
Finally, we apply the new 8th-order integrator to the one-dimensional Schrödinger equation (with $\hbar = 1$):
  \begin{equation*}
    i \frac{\partial}{\partial t}\phi(x,t)=-\frac{1}{2}\frac{\partial^2}{\partial x^2}\phi(x,t) + V(x)\phi(x,t),
  \end{equation*}
  where $\phi(x,t)$ is the wave function describing the physical system, and the potential is given by: 
  \begin{equation*}
    V(x)=-\frac{1}{2}x^2+\frac{1}{20}x^4.
  \end{equation*}
  Discretizing the spatial variable using a Fourier spectral collocation method yields an $N$-dimensional linear differential equation:
  \begin{equation*}
    i\dot{u}(t)= Hu(t) = (T+V)u(t), \qquad \text{with } u(t) \in \mathbb{C}^N,
  \end{equation*}
  where $T$ is the differentiation matrix resulting from the discretization of the kinetic energy operator,
  and $V$ is a diagonal matrix representing the potential, with $V_{nn} = V(x_n)$.
  The vector $u$ approximates the wave function: $u_n(t) \approx \phi(x_n, t)$.
The action of the matrix $T$ on $u$ is implemented using forward and backward discrete Fourier transforms \cite{lubich2008fqt}.

\begin{figure}[!h]
  \begin{center}
    \includegraphics[width=8cm]{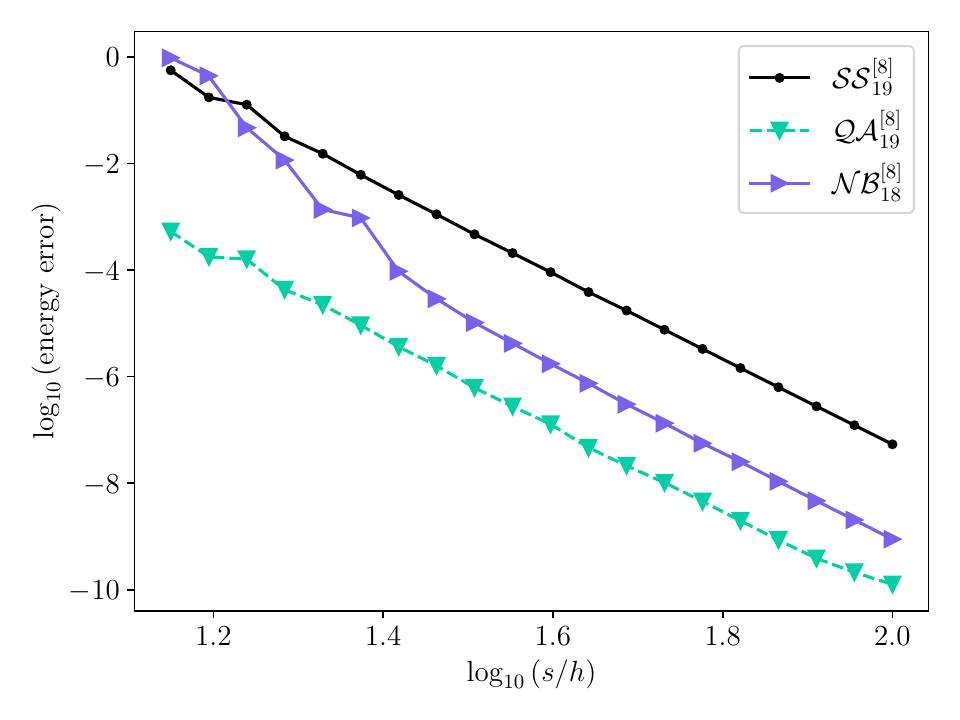}
    \caption{Efficiency diagram for the Schrödinger equation with a quartic potential using 8th-order methods. The final integration time is $t_f = 1000$.}
    \label{fig:5}
  \end{center}
\end{figure}

Figure~\ref{fig:5} shows an efficiency diagram for this test case, comparing our method $\mathcal{QA}_{19}^{[8]}$
against the most efficient known 8th-order RKN method and the best $\mathcal{SS}$-type method of the same order.
The initial condition is set as $u_n(0) = \sigma e^{-x_n^2 / 2}$, where $\sigma$ is a normalization constant.
The spatial domain is discretized over the interval $x \in [-10, 10]$ using $N = 256$ grid points, and time integration is performed up to $t_f = 1000$.
As a measure of accuracy, we compute the deviation in the expected energy at final time compared to the exact value:

  \begin{equation*}
    \text{energy error :}\quad|u^*(t_f)\cdot (Hu(t_f))- u^*(0)\cdot (Hu(0))|.
  \end{equation*}
  \section{Concluding remarks}
  In this work, we have shown, using the Lie formalism, why the number of order conditions can be significantly
  reduced for Hamiltonian problems with polynomial potentials. This theoretical insight has allowed us to construct
  new integrators that outperform the most efficient methods currently available in the literature for both cubic and quartic potentials.

  These new integrators should be regarded as highly competitive alternatives when the goal is to integrate
  cubic or quartic potentials while preserving the geometric structure of the original system.

  The methods presented here can be seen as a natural continuation of the work by \cite{iserles1998rkm},
  and they demonstrate that designing specialized integrators for structured problems can yield efficiency gains of several orders of magnitude.

  \section*{Acknowledgements}
  This work was supported by Ministerio de Ciencia e Innovación (Spain) under project PID2022- 136585NB-C21, MCIN/AEI/10.13039/501100011033/FEDER, UE and
  Universitat Jaume I under project 25I675.

  The author wish to thank to F. Casas and S. Blanes for their important comments and remarks.
%\bibliographystyle{abbrv}
%\bibliography{biblio}

\end{document}